\documentclass[twocolumn]{article}
\usepackage{amsmath,bm,url,mathrsfs,color}
\def\P{{\rm P}}        
\def\T{\textsf{T}}     
\allowdisplaybreaks

\title{A world record in Atlantic City\\ and the length of the shooter's hand at craps}
\author{S. N. Ethier\thanks{Department of Mathematics, University of Utah, 155 S. 1400 E., Salt Lake City, UT 84112, USA. e-mail: ethier@math.utah.edu}\; and Fred M. Hoppe\thanks{Department of Mathematics and Statistics, McMaster University, 1280 Main Street W., Hamilton, ON L8S 4K1, Canada.  e-mail: hoppe@mcmaster.ca}}
\date{}

\begin{document}
\maketitle

It was widely reported in the media that, on 23 May 2009, at the Borgata Hotel Casino \& Spa in Atlantic City, Patricia DeMauro\footnote{Spelled Demauro in some accounts.}, playing craps for only the second time, rolled the dice for four hours and 18 minutes, finally sevening out at the 154th roll.  Initial estimates of the probability of this event ranged from one chance in 3.5 billion \cite{Star-Ledger} to one chance in 1.56 trillion \cite{Time}.  Subsequent computations agreed on one chance in 5.6 (or 5.59) billion \cite{WSJblog,ChanceNews,Wizard}.

This established a new world record, previously held by the late Stanley Fujitake (118 rolls, 28 May 1989, California Hotel and Casino, Las Vegas) \cite{GoldenArm}.  One might ask how reliable these numbers (118 and 154) are.  In Mr.~Fujitake's case, casino personnel replayed the surveillance videotape to confirm the number of rolls and the duration of time (three hours and six minutes).   We imagine that the same happened in Ms.~DeMauro's case.

There is also a report that Mr.~Fujitake's record was broken earlier by a gentleman known only as The Captain (148 rolls, July 2005, Atlantic City) \cite[Part~4]{Scoblete}.  However, this incident is not well documented and was unknown to Borgata officials.  In fact, a statistical argument has been offered \cite[p.~480]{Grosjean} suggesting that the story is apocryphal.

Our aim in this article is not simply to derive a more accurate probability, but to show that this apparently prosaic problem involves some interesting mathematics, including Markov chains, matrix theory, and Galois theory.

\section*{Background}

Craps is played by rolling a pair of dice repeatedly.  For most bets, only the sum of the numbers appearing on the two dice matters, and this sum has distribution
\begin{equation}\label{dice}
\pi_j:={6-|j-7|\over36}, \qquad j=2,3,\ldots,12.
\end{equation}
The basic bet at craps is the \textit{pass-line bet}, which is defined as follows.  The first roll is the \textit{come-out roll}.  If 7 or 11 appears (a \textit{natural}), the bettor wins.  If 2, 3, or 12 appears (a \textit{craps number}), the bettor loses.  If a number belonging to 
$$
\mathscr{P}:=\{4,5,6,8,9,10\}
$$
appears, that number becomes the \textit{point}.  The dice continue to be rolled until the point is repeated (or \textit{made}), in which case the bettor wins, or 7 appears, in which case the bettor loses.  The latter event is called a \textit{seven out}.  A win pays even money.  The first roll following a decision is a new come-out roll, beginning the process again.

A shooter is permitted to roll the dice until he or she sevens out.  The sequence of rolls by the shooter is called the \textit{shooter's hand}.  Notice that the shooter's hand can contain winning 7s and losing decisions prior to the seven out.  The \textit{length} of the shooter's hand (i.e., the number of rolls) is a random variable we will denote by $L$.  Our concern here is with 
\begin{equation}\label{tail}
t(n):=\P(L\ge n),\qquad n\ge1,
\end{equation}
the tail of the distribution of $L$.   For example, $t(154)$ is the probability of achieving a hand at least as long as that of Ms.~DeMauro.  As can be easily verified from 
(\ref{recursion}), (\ref{matrixformula}), or (\ref{closed}) below, $t(154)\approx0.178\,882\,426\times10^{-9}$; 
to state it in the way preferred by the media, this amounts to one chance in 5.59 billion, approximately.  The 1 in 3.5 billion figure came from a simulation that was not long enough.  The 1 in 1.56 trillion figure came from $(1-\pi_7)^{154}$, which is the right answer to the wrong question.

\section*{Two methods}

We know of two methods for evaluating the tail probabilities (\ref{tail}).
The first is by recursion.  As pointed out in \cite{Ethier}, $t(1)=t(2)=1$ and 
\begin{eqnarray}\label{recursion}
\lefteqn{t(n)}\nonumber\\
&=&\bigg(1-\sum_{j\in\mathscr{P}}\pi_j\bigg)t(n-1)+\sum_{j\in\mathscr{P}}\pi_j(1-\pi_j-\pi_7)^{n-2}\nonumber\\
&&\quad{}+\sum_{j\in\mathscr{P}}\pi_j\sum_{l=2}^{n-1}(1-\pi_j-\pi_7)^{l-2} \pi_j
\,t(n-l)
\end{eqnarray}
for each $n\ge3$.  Indeed, for the event that the shooter
sevens out in no fewer than $n$ rolls to occur, consider the result of the initial
come-out roll.  If a natural or a craps number occurs, then, beginning with the
next roll, the shooter must seven out in no fewer than $n-1$ rolls.  If a
point number occurs, then there are two possibilities.  Either the point is still
unresolved after $n-2$ additional rolls, or it is made at roll $l$ for some
$l\in\{2,3,\ldots,n-1\}$ and the shooter subsequently sevens out in no fewer than
$n-l$ rolls. 

The second method, first suggested, to the best of our knowledge, by Peter A. Griffin in 1987 (unpublished) and rediscovered several times since, is based on a Markov chain.  The state space is 
\begin{equation}\label{S}
S:=\{{\rm co},{\rm p}4\mbox{-}10,{\rm p}5\mbox{-}9,{\rm p}6\mbox{-}8,7{\rm o}\} \equiv \{1,2,3,4,5\},
\end{equation}
whose five states represent the events that the shooter is coming out, has established the point 4 or 10, has established the point 5 or 9, has established the point 6 or 8, and has sevened out.  The one-step transition matrix, which can be inferred from (\ref{dice}), is
\begin{equation}\label{P}
\bm P:={1\over36}\left(\begin{array}{ccccc}
12&6&8&10&0\\
3&27&0&0&6\\
4&0&26&0&6\\
5&0&0&25&6\\
0&0&0&0&36\end{array}\right).
\end{equation}
The probability of sevening out in $n-1$ rolls or less is then just the probability that absorption in state 7o occurs by the $(n-1)$th step of the Markov chain, starting in state co.  A marginal simplification results by considering the 4 by 4 principal submatrix $\bm Q$ of (\ref{P}) corresponding to the transient states.
Thus, we have
\begin{equation}\label{matrixformula}
t(n)=1-(\bm P^{n-1})_{1,5} =  \sum_{j=1}^4(\bm Q^{n-1})_{1,j}.
\end{equation}
Clearly, (\ref{recursion}) is not a closed-form expression, and we do not regard (\ref{matrixformula}) as being in closed form either.  Is there a closed-form expression for $t(n)$?

\section*{Positivity of the eigenvalues}

We begin by showing that the eigenvalues of $\bm Q$ are positive.  The determinant of
\begin{eqnarray*}
\lefteqn{\bm Q - z \bm I}\\ &=& {1\over36}
\left(\begin{array}{cccc}
12- 36z&6&8&10\\
3&27- 36z&0&0\\
4&0&26- 36z&0\\
5&0&0&25- 36z\\
\end{array}\right).
\end{eqnarray*}
is unaltered by row operations. From the first row subtract $6/(27-36z)$ times the second row, $8/(26-36z)$ times the third row,
and $10/(25- 36z)$ times the fourth row,  cancelling the entries 6/36, 8/36, and 10/36 and making the (1,1) entry equal to 1/36 times
\begin{equation}\label{1,1-entry}
12- 36z -  3\, \frac{6}{ 27-36z} - 4\, \frac{8}{26-36z} - 5\, \frac{10}{25- 36z}.
\end{equation}
The determinant of $\bm Q - z \bm I $, and therefore the characteristic polynomial $q(z)$ of 
$\bm Q$ is then just the product of the diagonal entries in the transformed matrix,
which is (\ref{1,1-entry}) multiplied by $(27- 36z)(26- 36z)(25- 36z)/(36)^4$.  Thus,
\begin{eqnarray*}
q(z)&=&[(12- 36z)(27- 36z)(26- 36z)(25- 36z)\\
&&\quad{}-18(26- 36z)(25- 36z)\\
&&\quad{}-32(27- 36z)(25- 36z)\\
&&\quad{}-50(27- 36z)(26- 36z)]/(36)^4.
\end{eqnarray*}
We find that $q(1),q(27/36),q(26/36),q(25/36),q(0)$ alternate signs and therefore the eigenvalues are positive and interlaced between the diagonal entries (ignoring the entry 12/36).  
More precisely, denoting the eigenvalues by $1>e_1>e_2>e_3>e_4>0$, we have
$$
1>e_1>{27\over36}>e_2>{26\over36}>e_3>{25\over36}>e_4>0.
$$

The matrix $\bm Q$, which has the structure of an arrowhead matrix \cite{HornJohnson}, is positive definite, although not symmetric.
This is easily seen by applying the same type of row operations to the symmetric part $\bm A = \frac{1}{2}(\bm Q+\bm Q^\T)$ of $\bm Q$ to show that the eigenvalues of $\bm A$ interlace its diagonal elements (except 12/36).  But a symmetric matrix is positive definite if and only if all its eigenvalues are positive, and a non-symmetric matrix is positive definite if and only if its symmetric part is positive definite, confirming that $\bm Q$ is positive definite.

\section*{A closed-form expression}

The eigenvalues of $\bm Q$ are the four roots of the quartic equation $q(z)=0$ or
$$
23328z^4-58320z^3+51534z^2-18321z+1975=0,
$$
while $\bm P$ has an additional eigenvalue, 1, the spectral radius.  
We can use the quartic formula (or \textit{Mathematica}) to find these roots.  We notice that the complex number 
$$
\alpha:=\zeta^{1/3}+{9829\over \zeta^{1/3}}, 
$$
where
$$
\zeta:=-710369+18i\sqrt{1373296647},
$$ 
appears three times in each root.  Fortunately, $\alpha$ is positive, as we see by writing $\zeta$ in polar form, that is, $\zeta=re^{i\theta}$.  We obtain
$$
\alpha=2\sqrt{9829}\,\cos\bigg[{1\over3}\cos^{-1}\bigg(-{710369\over9829\sqrt{9829}}\bigg)\bigg].
$$
The four eigenvalues of $\bm Q$ can be expressed as
\begin{eqnarray*}
e_1&:=&e(1,1),\\
e_2&:=&e(1,-1),\\
e_3&:=&e(-1,1),\\
e_4&:=&e(-1,-1),
\end{eqnarray*}
where
\begin{eqnarray*}
e(u,v)&:=&{5\over8}+{u\over72}\sqrt{{349+\alpha\over3}}\\
&&\quad{}+{v\over72}\sqrt{{698-\alpha\over3}-2136u\sqrt{3\over349+\alpha}}.
\end{eqnarray*}

Next we need to find right eigenvectors corresponding to the five eigenvalues of $\bm P$.  Fortunately, these eigenvectors can be expressed in terms of the eigenvalues.  Indeed,
with $\bm r(x)$ defined to be the vector-valued function
$$
\left(\begin{array}{c}-5+(1/5)x\\
-175+(581/15)x-(21/10)x^2+(1/30)x^3\\
275/2-(1199/40)x+(8/5)x^2-(1/40)x^3\\
1\\
0\end{array}\right)
$$
we find that right eigenvectors corresponding to eigenvalues $1,e_1,e_2,e_3,e_4$ are
$$
(1,1,1,1,1)^\T,\;\bm r(36e_1),\;\bm r(36e_2),\;\bm r(36e_3),\;\bm r(36e_4),
$$
respectively.  Letting $\bm R$ denote the matrix whose columns are these right eigenvectors and putting $\bm L:=\bm R^{-1}$, the rows of which are left eigenvectors, we know by (\ref{matrixformula}) and the spectral representation that
$$
t(n)=1-\{\bm R\,{\rm diag}(1,e_1^{n-1},e_2^{n-1},e_3^{n-1},e_4^{n-1})\bm L\}_{1,5}.
$$

After much algebra (and with some help from \textit{Mathematica}), we obtain
\begin{equation}\label{closed}
t(n)=c_1 e_1^{n-1}+c_2 e_2^{n-1}+c_3 e_3^{n-1}+c_4 e_4^{n-1},
\end{equation}
where the coefficients are defined in terms of the eigenvalues and the function
\begin{eqnarray*}
f(w,x,y,z)&:=&(-25+36w)[4835-5580(x+y+z)\\
&&\quad{}+6480(xy+xz+yz)-7776xyz]\\
&&\;\qquad/[38880(w-x)(w-y)(w-z)]
\end{eqnarray*}
as follows:
\begin{eqnarray*}
c_1&:=&f(e_1,e_2,e_3,e_4),\\ c_2&:=&f(e_2,e_3,e_4,e_1),\\ c_3&:=&f(e_3,e_4,e_1,e_2),\\ c_4&:=&f(e_4,e_1,e_2,e_3).
\end{eqnarray*}
Of course, (\ref{closed}) is our closed-form expression.  

Incidentally, the fact that
$t(1)=t(2)=1$ implies that
\begin{equation}\label{sum1}
c_1+c_2+c_3+c_4=1
\end{equation}
and
$$
c_1e_1+c_2e_2+c_3e_3+c_4e_4=1.
$$

In a sequence of independent Bernoulli trials, each with success probability $p$, the number of trials $X$ needed to achieve the first success has the \textit{geometric distribution} with parameter $p$, and
$$
P(X\ge n)=(1-p)^{n-1},\qquad n\ge1.
$$
It follows that the distribution of $L$ is \textit{a linear combination of four geometric distributions}.  It is not a convex combination: (\ref{sum1}) holds but, as we will see,
$$ 
c_1>0,\quad  c_2<0,\quad  c_3<0,\quad  c_4<0.  
$$
In particular, we have the inequality
\begin{equation}\label{ineq}
t(n)<c_1 e_1^{n-1},\qquad n\ge1,
\end{equation}
as well as the asymptotic formula
\begin{equation}\label{asymp}
t(n)\sim c_1 e_1^{n-1}\quad {\rm as\ }n\to\infty.
\end{equation}

\section*{Numerical approximations}

Rounding to 18 decimal places, the nonunit eigenvalues of $\bm P$ are
\begin{eqnarray*}
e_1&\approx&0.862\,473\,751\,659\,322\,030,\\
e_2&\approx&0.741\,708\,271\,459\,795\,977,\\
e_3&\approx&0.709\,206\,775\,794\,379\,015,\\
e_4&\approx&0.186\,611\,201\,086\,502\,979,
\end{eqnarray*}
and the coefficients in (\ref{closed}) are
\begin{eqnarray*}
c_1&\approx&\phantom{-}1.211\,844\,812\,464\,518\,572,\\
c_2&\approx&-0.006\,375\,542\,263\,784\,777,\\
c_3&\approx&-0.004\,042\,671\,248\,651\,503,\\
c_4&\approx&-0.201\,426\,598\,952\,082\,292.
\end{eqnarray*}
These numbers will give very accurate results over a wide range of values of $n$. 

The result (\ref{asymp}) shows that the leading term in (\ref{closed}) may be adequate for large $n$; it can be shown that
$$
1<c_1 e_1^{n-1}/t(n)<1+10^{-m}
$$
for $m=3$ if $n\ge19$; for $m=6$ if $n\ge59$; for $m=9$ if $n\ge104$; and for $m=12$ if $n\ge150$.

\section*{Crapless craps}

In crapless craps \cite[p.~354]{ScarneRawson}, as the name suggests, there are no craps numbers and 7 is the only natural.  Therefore, the set of possible point numbers is 
$$
\mathscr{P}_0:=\{2,3,4,5,6,8,9,10,11,12\}
$$
but otherwise the rules of craps apply.  

With $L_0$ denoting the length of the shooter's hand, the analogues of (\ref{S})--(\ref{matrixformula}) are
\begin{eqnarray*}
S_0&:=&\{{\rm co},{\rm p}2\mbox{-}12,{\rm p}3\mbox{-}11,{\rm p}4\mbox{-}10,{\rm p}5\mbox{-}9,{\rm p}6\mbox{-}8,7{\rm o}\}\\
&\;\equiv&\{1,2,3,4,5,6,7\},
\end{eqnarray*}
$$
\bm P_0:={1\over36}\left(\begin{array}{ccccccc}
6\;&2&4&6&8&10&0\\
1\;&29&0&0&0&0&6\\
2\;&0&28&0&0&0&6\\
3\;&0&0&27&0&0&6\\
4\;&0&0&0&26&0&6\\
5\;&0&0&0&0&25&6\\
0\;&0&0&0&0&0&36\end{array}\right),
$$
and 
$$
t_0(n):=\P(L_0\ge n)=1-(\bm P_0^{n-1})_{1,7}.
$$

There is an interesting distinction between this game and regular craps.  The nonunit eigenvalues of $\bm P_0$ are the roots of the sextic equation
\begin{eqnarray*}
&&15116544z^6-59206464z^5+93137040z^4\\
&&\;{}-73915740z^3+30008394z^2-5305446z+172975\\
&&\qquad{}=0,
\end{eqnarray*}
and the corresponding Galois group is, according to \textit{Maple}, the symmetric group $S_6$.  This means that our sextic is not solvable by radicals.  Thus, it appears that there is no closed-form expression for $t_0(n)$.

Nevertheless, the analogue of (\ref{closed}) holds (with six terms).  All nonunit eigenvalues belong to $(0,1)$ and all coefficients except the leading one are negative.  Thus, the analogues of (\ref{ineq}) and (\ref{asymp}) hold as well.  Also, the distribution of $L_0$ is \textit{a linear combination of six geometric distributions}.  These results are left as exercises for the interested reader.

Finally, $t_0(154)\approx0.296\,360\,068\times 10^{-10}$, which is to say that a hand of length 154 or more is only about one-sixth as likely as at regular craps (one chance in 33.7 billion, approximately).

\section*{Acknowledgment}

We thank Roger Horn for pointing out the interlacing property of the eigenvalues of $\bm Q$.

\bigskip



\end{document}